\documentclass[jsco]{academic}

\usepackage{amsmath}
\usepackage{amsfonts}
\usepackage{amscd}
\usepackage{diagrams}
\usepackage{epsfig}

\newdefinition{example}[theorem]{Example}

\begin{document}

\shortauthor{A. Khetan}
\shorttitle{The Resultant of an Unmixed Bivariate-System}
\title{The Resultant of an Unmixed Bivariate~System}
\author{Amit Khetan \footnote{E-mail: akhetan@math.berkeley.edu}}
\address{Department of Mathematics, UC Berkeley, Berkeley CA, USA}
\maketitle

\begin{abstract}

This paper gives an explicit method for computing the resultant of any
sparse unmixed bivariate system with given support. We construct
square matrices whose determinant is exactly the resultant. The
matrices constructed are of hybrid Sylvester and B\'ezout type. The
results extend those in \citep{Khe} by giving a complete combinatorial
description of the matrix. Previous work by D'Andrea \citep{D} gave
pure Sylvester type matrices (in any dimension).  In the bivariate
case, D'Andrea and Emiris \citep{DE} constructed hybrid matrices with
one B\'ezout row. These matrices are only guaranteed to have
determinant some multiple of the resultant. The main contribution of
this paper is the addition of new B\'ezout terms allowing us to
achieve exact formulas. We make use of the exterior algebra techniques
of Eisenbud, Fl{\o}ystad, and Schreyer \citep{ES, EFS}.

\end{abstract}

\section{Introduction}

Let $f_1, \dots, f_{n+1} \in \mathbb{C}[x_1, x_1^{-1}, \dots, x_n,
x_n^{-1}]$ be Laurent polynomials in $n$ variables with the same
Newton polytope $Q \subset \mathbb{R}^n$. Let $A = Q \cap
\mathbb{Z}^n$. So we can write:

\[ f_i = \sum_{\alpha \in A} C_{i\alpha}x^{\alpha}. \]

\noindent We will assume that $Q$ is actually $n$-dimensional, and furthermore
that $A$ affinely spans $\mathbb{Z}^n$.

\begin{definition}
\label{d:ares}
The $A$-{\em resultant} ${\rm Res}_A(f_1, \dots, f_{n+1})$ is the
irreducible polynomial in the $C_{i\alpha}$, unique up to sign, which
vanishes whenever $f_1, \dots, f_{n+1}$ have a common root in the
algebraic torus $(\mathbb{C}^\ast)^n$.
\end{definition}

The existence, uniqueness, and irreducibility of the $A$-resultant are
proved in the book by Gelfand, Kapranov, and Zelevinsky \citep{GKZ}.
The $A$-resultant, also called the sparse resultant, allows one to
eliminate $n$ variables from $n+1$ unmixed equations. Hence,
resultants can be quite useful in solving systems of polynomial
equations \citep{CLO2}. It is an important problem to find
efficiently computable, explicit formulas for the resultant.

When $n=1$, we are in the case of the classical resultant of two
polynomials in one variable of the same degree. There are two formulas
due to Sylvester and B\'ezout which represent the resultant as the
determinant of an easily computable matrix. Sylvester's matrix has
entries that are either 0 or a coefficient of $f_1$ or $f_2$. The
entries in B\'ezout's matrix are linear in the coefficients of each of
the $f_i$ hence quadratic overall. 

Our work deals with the case $n=2$. We give a determinantal formula
which is of hybrid Sylvester and B\'ezout type. A preliminary version of
these results appeared in the ISSAC 2002 Proceedings \citep{Khe}. This paper
makes the formula completely explicit and provides complete proofs.
Our approach follows work by Jouanolou \citep{J} and Dickenstein and
D'Andrea \citep{DD} who found formulas for the ``dense'' resultant,
when the polytope $Q$ is a coordinate simplex of some degree. We make
heavy use of new techniques by Eisenbud, Fl{\o}ystad and Schreyer
\citep{ES, EFS} relating resultants to complexes over an exterior algebra.

\begin{theorem}
\label{thm:blmtrx}
The resultant of a system $(f_1, f_2, f_3) \in \mathbb{C}[x_1, x_2,
x_1^{-1}, x_2^{-1}]$ with common Newton polygon $Q$ is the determinant
of the block matrix:
\begin{equation*}
\label{e:blmtrx}
\begin{pmatrix}
 B & L \\ \tilde{L} & 0 \end{pmatrix}.
\end{equation*}
\noindent The entries of $L$ and $\tilde{L}$ are linear forms, and the entries
of $B$ are cubic forms in the coefficients
$C_{i\alpha}$.
\end{theorem} 

The columns of $B$ and $\tilde{L}$ are indexed by the lattice points in
$Q$, the rows of $B$ and $L$ are indexed by the interior
lattice points in $2 \cdot Q$, the matrix $\tilde{L}$ has three rows indexed
by $\{ f_1, f_2, f_3 \}$, and the columns of the matrix $L$ are
indexed by pairs $(f_i, a)$ where $i \in \{1,2,3\}$ and $a$ runs over
the interior lattice points of $Q$. 
Each entry of $L$ and $\tilde L$ is either zero or is a coefficient of
some $f_i$ and is determined in the following straightforward
manner. The entry of $\tilde L$ in row $f_i$ and column $a$ is the
coefficient of $x^a$ in $f_i$. The entry of $L$ in row $b$ and
column $(f_i, a)$ is the coefficient of $x^{b-a}$ in $f_i$. The
entries of the matrix $B$ are linear forms in {\em bracket variables}.
A bracket variable is defined as

$$[abc] = \det \bmatrix  C_{1a} & C_{1b} & C_{1c} \\ 
                         C_{2a} & C_{2b} & C_{2c} \\
                         C_{3a} & C_{3b} & C_{3c} \endbmatrix,$$

\noindent where $C_{ia}$ is the coefficient of $x^a$ in $f_i$. An explicit
formula for $B$ is described in Section 3 below.

\begin{figure}
\begin{center}
\epsfig{file=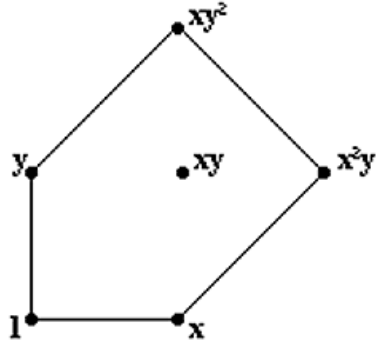}
\caption{Newton Polygon Q}
\label{f:newt}
\end{center}
\end{figure}

\begin{example}

\label{example}

\begin{align*}
f_1 &= C_{11} + C_{12}x + C_{13}y + C_{14}xy + C_{15}x^2y + C_{16}xy^2 \\
f_2 &= C_{21} + C_{22}x + C_{23}y + C_{24}xy + C_{25}x^2y + C_{26}xy^2\\
f_3 &= C_{31} + C_{32}x + C_{33}y + C_{34}xy + C_{35}x^2y + C_{36}xy^2
\end{align*}

The system above has Newton polygon as shown in Figure
\ref{f:newt}. We will show that the resultant of this system is the
determinant of the matrix in Table \ref{tbl:matrix}.

\begin{table}
\begin{center}
\begin{tabular}{|c c c c c c c c c|}
0	& $[124]$       &     0     & $[126]-[234]$ & $-[235]$	    & $-[236]$ & $c_{11}$ & $c_{21}$ & $c_{31}$ \\
0       &    0          &     0     &   0           &   0           &	   0   & $c_{12}$ & $c_{22}$ & $c_{32}$ \\
0	& $[126]-[135]$ &     0     & $[146]-[236]$ & $[156]+[345]$ & $[346]$  & $c_{13}$ & $c_{23}$ & $c_{33}$ \\
0       & $-[145]$      &     0     & $[156]-[345]$ & $[256]$       & $[356]$  & $c_{14}$ & $c_{24}$ & $c_{34}$ \\
0       &    0          &     0     &    0          &	 0          &    0     & $c_{15}$ & $c_{25}$ & $c_{35}$ \\
0       & $[156]$       &     0     & $[356]$       & $[456]$	    &	 0     & $c_{16}$ & $c_{26}$ & $c_{36}$ \\
$c_{11}$& $c_{12}$      &  $c_{13}$ & $c_{14}$      & $c_{15}$      & $c_{16}$ &    0	  &    0     &    0     \\
$c_{21}$& $c_{22}$      &  $c_{23}$ & $c_{24}$      & $c_{25}$      & $c_{26}$ &    0	  &    0     &    0     \\
$c_{31}$& $c_{32}$      &  $c_{33}$ & $c_{34}$      & $c_{35}$      & $c_{36}$ &    0	  &    0     &    0     \\
\end{tabular}
\end{center}
\caption{Resultant matrix for Example \ref{example}}
\label{tbl:matrix}
\end{table}

\end{example}

In Section \ref{s:toric} we provide some preliminary results about
toric varieties and their homogeneous coordinates which allow us to
present our formula in Section \ref{s:formula}. Section \ref{s:tate}
describes the exterior algebra techniques of Eisenbud, Schreyer, and
Fl{\o}ystad. Section \ref{s:torictate} applies these results to the
toric setting, while Section \ref{s:proofs} goes on to prove our
formula.  Finally Section \ref{s:general} briefly discusses possible
generalizations to more variables.

\subsection{Acknowledgments}

I would like to thank my advisor Bernd Sturmfels for providing
direction and support. I thank David Eisenbud for introducing me to
exterior algebra methods. I thank David Speyer for the proof of Lemma
\ref{comblemma}.

\section{Toric Varieties}
\label{s:toric}

\begin{definition} 
Let $\mathbb{Q} \subset \mathbb{R}^n$ be a lattice polytope of
dimension $n$, and $A = Q \cap \mathbb{Z}^n = \{ \alpha_1, \dots,
\alpha_N \}$. We assume that $A$ affinely spans $\mathbb{Z}^n$. The
{\em toric variety} $X_A$ is the dimension $n$ variety defined as the
Zariski closure of the following set in $\mathbb{P}^{N-1}$:

\[ X_A^0 = \{ (x^{\alpha_1}:\cdots : x^{\alpha_N}) \ : \ x = (x_1, \dots, x_n) \in (\mathbb{C}^\ast)^n \}. \]

\end{definition}

Now a polynomial system $(f_1, \dots, f_{n+1})$ can be thought of as
$n+1$ hyperplane sections of $X_A$ in $\mathbb{P}^{N-1}$. Generically,
such a system defines a codimension $n+1$ plane.

For any $n$-dimensional irreducible projective variety $X$, it turns
out that the condition on a linear subspace of codimension $n+1$
meeting $X$ is actually a closed condition of codimension 1 (see
\citep{GKZ} for details). Therefore we can make the following
definition.

\begin{definition} 
If $X \subset \mathbb{P}^{N-1}$ is a variety of dimension $n$, the
codimension $n+1$ planes meeting $X$ define a hypersurface in the
Grassmannian $G(n+1, N)$. The equation of this hypersurface is called
the {\em Chow form} of $X$.
\end{definition}

In particular, the $A$-resultant is the Chow form of $X_A$.  As a
consequence we have the following strengthening of Definition
\ref{d:ares}.

\begin{corollary} 
${\rm Res}_A(f_1, \dots, f_{n+1}) = 0$ if and only if the $f_i$ have a
common root on $X_A$.
\end{corollary} 

Returning to the defining polytope $Q$, let $d_1, \dots, d_s$ denote
the facets of $Q$. Let $\eta_i$ be the first lattice vector along the
inner normal to facet $d_i$. The {\em normal fan} of $Q$ is the set of
cones, one for each vertex, spanned by the $\eta_i$ corresponding to
facets incident to that vertex.  The next proposition can be found in
Fulton's book \citep{Ful}. 

\begin{proposition} 
The $\eta_i$ are in 1-1 correspondence with the $T$-invariant prime Weil
 divisors on $X_A$. Let $D_i$ denote the divisor corresponding to
 $\eta_i$.
\end{proposition}

The polytope $Q$ can be characterized completely in terms of the rays
in its normal fan as follows:

\[ Q = \{ m \in \mathbb{R}^n \: \ \langle m, \eta_i \rangle \geq -a_i, \ i = 1, \dots, s \}. \]

 The very ample divisor corresponding to the embedding of $X_A$ into
$\mathbb{P}^{N-1}$ corresponding to $Q$ is just $D = \sum a_i D_i$.
We can now define the {\em homogeneous coordinate ring} of $X_A$. This
was introduced by Cox \citep{Cox} and the propositions below follow
from this paper.

 Let $S = \mathbb{C}[y_1, \dots, y_s]$ be the polynomial ring with one
variable for each $\eta_i$. Consider the short exact sequence of
abelian groups:

$$\begin{CD} 0 @>>> \mathbb{Z}^{n} @>\phi>> \mathbb{Z}^s @>\pi>> G
@>>> 0. \end{CD}$$

\noindent Here $\phi(m) = (\langle m, \eta_1 \rangle, \dots, \langle m, \eta_s
\rangle)$, and $G$ is the cokernel of $\phi$.

\begin{definition}
Define a $G$-grading on $S$ as follows. Given $y^\alpha \in S$, let
$\deg(y^\alpha) = \pi(\alpha) \in G$.
\end{definition}

Now we will identify the lattice points in $Q$ with a graded piece of
$S$.

\begin{definition} 
Let $\alpha \in Q \cap \mathbb{Z}^n$. Define $\alpha_i = \langle
\alpha, \eta_i \rangle + a_i$ for $i = 1, \dots, s$ and the $a_i$ are
the defining data for $Q$ as above. The $Q$-{\em homogenization} of
$x^\alpha$ is $\prod_{i=1}^s y_i^{\alpha_i}$. We will write this as
$y^{\alpha}$ and use the letter $\alpha$ to denote both a vector
$\alpha \in \mathbb{Z}^n$ and its homogenization $(\alpha_1, \dots
\alpha_s)$, where the meaning will be clear from the context.
\end{definition}

\begin{proposition} 
\label{p:sq}
Let $a = (a_1, \dots, a_s)$ be the defining data for $Q$. The
monomials in the $\pi(a)$ graded piece of $S$ are in 1-1
correspondence with the lattice points in $Q$. Denote this graded
piece by $S_Q$. Moreover, $H^0(X_A, \mathcal{O}(D)) \cong S_Q$.
\end{proposition}

There is a similar characterization of the interior lattice points of
$Q$.

\begin{proposition} 
\label{p:sintq}
Let $\omega_0 = (1, 1, \dots, 1) \in \mathbb{Z}^s$. The monomials in the
$\pi(a - \omega_0)$ graded piece of $S$ are in 1-1 correspondence with
the interior lattice points of $Q$. Denote this graded piece $S_{{\rm
int}(Q)}$. We have $H^0(X, \mathcal{O}(D- \sum_{i=1}^s D_i)) \cong
S_{{\rm int}(Q)}$.
\end{proposition}

\section{Formula for $B$}
\label{s:formula}

We now return to case of two variables. So $(f_1, f_2, f_3) \in
\mathbb{C}[x_1, x_2, x_1^{-1}, x_2^{-1}]$ have common Newton polygon
$Q \subset \mathbb{R}^2$. The rays in the normal fan of $Q$ are $\{
\eta_1, \dots, \eta_s \}$, assumed to be in counterclockwise order.  We
pick out the distinguished cone spanned by $\{\eta_1, \eta_2\}$ and
partition the vectors in the fan as follows:

\begin{align} 
\label{eqn:Rsets}
R_1 &= \{i \ | \ \eta_i = c_1 \eta_1 + c_2 \eta_2 \ {\rm with} \  c_1 \geq 0 \ {\rm and} \ c_2 \leq 0 \} \nonumber \\
R_2 &= \{i \ | \ \eta_i = c_1 \eta_1 + c_2 \eta_2 \ {\rm with} \  c_1 \leq 0 \ {\rm and} \ c_2 \geq 0 \} \\
R_3 &= \{i \ | \ \eta_i = c_1 \eta_1 + c_2 \eta_2 \ {\rm with} \  c_1 < 0 \ {\rm and} \ c_2 < 0 \}. \nonumber 
\end{align}

It is possible that $R_3$ as defined is empty. If that is the case we
need to {\em refine} the fan, by adding in one new vector, say
$\eta_{s+1} = -\eta_1 - \eta_2$. This new vector $\eta_{s+1}$ lies in
the interior of some cone spanned by $\eta_i$ and $\eta_{j}$, hence
can be written as $c_1 \eta_i + c_2 \eta_j$ for some positive $c_1,
c_2$. Define $a_{s+1} = c_1 a_i + c_2 a_j$.  As above, given $\alpha
\in Q$ we denote by $\alpha_{s+1}$ the quantity $\langle \alpha,
\eta_{s+1} \rangle + a_{s+1}$.

In fact, if there is a single fan vector $\eta_i$
such that $-\eta_i$ is not a ray in the fan, then we can choose our
distinguished cone to be the one containing $-\eta_i$, and $R_3$ is
guaranteed not to be empty. However, for polytopes such that every
edge has a corresponding parallel edge, this is not the case.

A good way to think about these sets is that we choose a distinguished
vertex $p$ of $Q$ having normal cone spanned by $\{\eta_1,
\eta_2\}$. The set $R_3$ consists of all edges of $Q$ such that the
corresponding inner normals are maximized at $v$. If there is no such
edge, then our refinement adds in a ``length 0'' edge whose inner
normal is maximized at $p$. $R_1$ is the set of the remaining edges
clockwise from $v$, while $R_2$ is the set of remaining edges
counterclockwise from $v$.

\begin{figure}
\epsfig{file=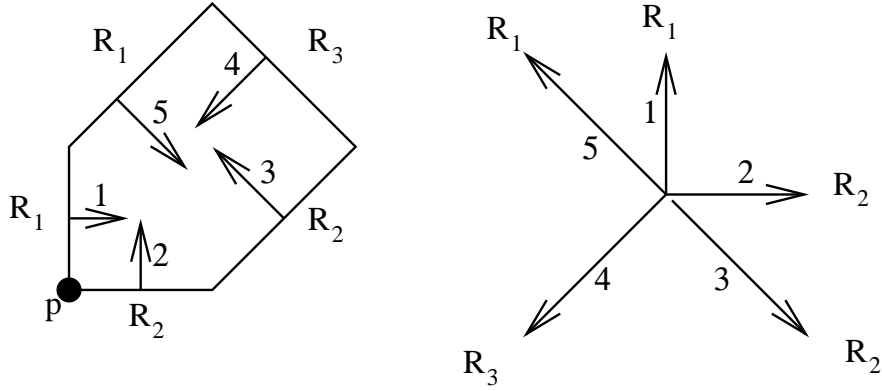}
\caption{The normal fan and partition $R_1, R_2, R_3$}
\label{f:fan}
\end{figure}

This partition is illustrated in Figure \ref{f:fan} for Example \ref{example}
with the choice of the vertex $p$.  Edge $4$ has the only
normal maximized at $p$, thus is the only element in $R_3$. The edges
in $R_1$ and $R_2$ are $\{1,5\}$ and $\{2, 3\}$ respectively. 

We can now state an explicit formula for the matrix $B$ appearing in
the Theorem \ref{thm:blmtrx}.

\begin{theorem}
\label{thm:main}

 The matrix $B$ from Theorem \ref{thm:blmtrx} is the matrix of the linear map
 $\Delta_Q \ : \ (S_Q)^\ast \to S_{{\rm int}(2Q)}$ defined as follows:

 \[ \Delta_Q((y^\alpha)^\ast) = \sum_{(u,v,w) \in F_{\alpha} \subset
 A^3} [uvw] y^{u+v+w - \alpha - \omega_0}. \]

Here $\omega_0 = (1, 1, \dots, 1)$, and $F_\alpha$ is the set of all triples
$(u,v,w) \in A^3$ satisfying the following Boolean combination of inequalities:
\begin{align}
\label{eqn:ineq}
\forall i \in R_1 \quad u_i + v_i + w_i >& \alpha_i \nonumber \\
\exists i \in R_1 \quad v_i + w_i \leq& \alpha_i \nonumber \\
\forall j \in R_2 \quad v_j + w_j >& \alpha_j \\
\exists j \in R_2 \quad w_j \leq& \alpha_j \nonumber \\
\forall k \in R_3 \quad w_k >& \alpha_k \nonumber,
\end{align}
\noindent where the $R_i$ are as described in (\ref{eqn:Rsets}).

\end{theorem}

\begin{example}
Let's see how this works for Example \ref{example}. Specifically,
consider the point $\alpha = (1,1)$ corresponding to the monomial
$xy$.  The homogenization is $y_1y_2y_3y_4y_5$.  If the monomials are
numbered $1, \dots, 6$ as in the equations, then the only solutions to
the inequalities above are:
$$ (u,v,w) = \{ (2,6,1), (4, 6, 1), (5, 6, 1), (2,4,3), (5,4,3), (2,6,3),
(5,6,3) \}.$$

\noindent It follows that
\begin{align*} \Delta_Q((y_1y_2y_3y_4y_5)^\ast) &= ([261]+[243])y_3y_4^3y_5 + ([461] + [263])y_2y_3^2y_4^2 \\
&+ ([561] + [543])y_1y_2y_3y_4y_5 + [563]y_1y_2^2y_3^2, 
\end{align*}
\noindent which corresponds to the fourth column of the matrix in Table \ref{tbl:matrix}.

\end{example}

\section{Tate Resolution}

\label{s:tate}
In this section we describe a complex used by Eisenbud and Schreyer
\citep{ES,EFS} to compute Chow forms of projective varieties. This
begins as a complex of free modules over an exterior algebra, however
there is a functor which transforms it into a complex of vector
bundles on the Grassmannian. The determinant of this new complex will
be the Chow form.

Suppose $X \subset \mathbb{P}^{N-1}$ is an irreducible variety of
dimension $ n$. The ambient projective space $\mathbf{P} =
\mathbb{P}^{N-1}$ has the graded coordinate ring $R = \mathbb{C}[X_1,
\dots, X_N]$. If we let $W$ be the $\mathbb{C}$ vector space spanned
by the $X_i$, identified with the degree 1 part of $R$, then
$\mathbf{P}$ is the projectivization $\mathbb{P}(W)$. The ring $R$ can
also be identified with the symmetric algebra $Sym(W)$.

Now let $V = W^\ast$, the dual vector space, with a corresponding dual
basis $ e_1, \dots e_N$.  We will consider the {\em exterior algebra}
$E = \bigwedge V$, which is also graded where the $e_i$ have degree $
-1$.  We will use the standard notation $E(k)$ to refer the rank
1 free $E$-module generated in degree $-k$.

Now given any coherent sheaf $\mathcal{F}$ on $\mathbf{P}$, there is
an associated exact complex of graded free $E$-modules, called the
{\em Tate resolution}, denoted $T(\mathcal{F})$. The terms of
$T(\mathcal{F})$ can be written in terms of the sheaf cohomology of
twists of $\mathcal{F}$. Namely, we have:

\begin{equation}
\label{e:tate}
T^e(\mathcal{F}) =  \oplus_{j=0}^{N-1} [H^j(\mathcal{F}(e-j))
\otimes_{\mathbb{C}} E(j-e)]
\end{equation}

\noindent for all $e \in \mathbb{Z}$. See
Eisenbud-Fl{\o}ystad-Schreyer \citep{EFS}.

Suppose further that $\mathcal{F}$ is chosen to be supported on
$X$. Recall that the Chow form of $X$ is the defining equation of the
set of codimension $n+1$-planes meeting $X$. Such a plane is specified
by an $n+1$ dimensional subspace $W_f = \mathbb{C} \{ f_1, \dots,
f_{n+1} \} \subset W$. Let $\mathbf{G_{n+1}}$ be the Grasmannian of
codimension $n+1$-planes on $\mathbf{P}$. Let $\mathcal{T}$ be the
{\em tautological bundle} on $\mathbf{G_{n+1}}$, that is to say the fiber at
the point corresponding to $f$ is just $W_f$.

The following proposition is a
consequence of Theorem 0.1 in \citep{EFS}.

\begin{proposition}

There is an additive functor $U_{n+1}$ from graded free modules over
$E$ to vector bundles on $\mathbf{G_{n+1}}$, such that $U_{n+1}(E(p))
= \wedge^p \mathcal{T}$. Furthermore, if $\mathcal{F}$ is a sheaf of
rank $k$ supported on a variety $X \subset \mathbb{P}(V)$ of dimension
$n$, $U_{n+1}(T(\mathcal{F}))$ is a complex of vector bundles whose
determinant is the $k$-th power of the Chow form of $X$.
\end{proposition}

The determinant of a complex of vector bundles on $\mathbf{G_{n+1}}$
is a homogeneous polynomial function on $\mathbf{G_{n+1}}$ whose value
at a particular point is the corresponding determinant of the complex
of {\em vector spaces} over that point. The determinant of a complex
of vector spaces is defined in \citep[Appendix A]{GKZ}. 

So, in particular if we could choose $\mathcal{F}$ so that enough
cohomology vanishes, this new complex $U_{n+1}(T(\mathcal{F}))$ may
have only two terms and a single non trivial map $\Psi_{\mathcal{F}}$.
Such sheaves are called {\em weakly Ulrich}, see \citep[Section
2]{ES}. In this case, to compute the Chow form we need only compute
the determinant of $\Psi_{\mathcal{F}}$. This is exactly what we do in
the next section. But first we need to describe the maps in the Tate
resolution, and also how the functor $U_{n+1}$ acts.

The maps in the Tate resolution are composed of maps
$H^j(\mathcal{F}(e-j)) \otimes E(j-e) \to H^k(\mathcal{F}(e+1-k))
\otimes E(k-e-1)$. All such maps for $k > j$ must be 0 by degree
considerations. 

When $k = j$ we have a {\em linear map} $H^j(\mathcal{F}(e-j)) \otimes
E(j-e) \to H^j(\mathcal{F}(e+1-j)) \otimes E(j-e-1)$ which is
canonical and completely well understood. Explicitly we consider the
graded $R$-module $M^j = \oplus_{l > 0} H^j(\mathcal{F}(l))$.  The
Bernstein-Gel'fand-Gel'fand correspondence \citep[Section
2]{EFS} applied to $M^j$ results in a map $M^j_{e-j} \otimes E(j-e)
\to M^j_{e-j+1} \otimes E(j-e-1)$ which is just multiplication by the
element $m = \sum X_i \otimes e_i$. By \citep[Theorem 4.1]{EFS} these
are exactly the linear maps in the Tate Resolution.

Much more mysterious are the nonlinear {\em diagonal maps}
corresponding to $k < j$. Indeed one of the major contributions of
this paper is an explicit formula for one of these diagonal maps in
the case of a toric surface. Eisenbud and Schreyer \citep{ES} outline a
general procedure for computing the Tate resolution, and therefore the
diagonal maps, however it requires computing a free resolution and is
not an explicit formulation.

Before moving on to the toric setting let us complete the description
of the functor $U_{n+1}$ by describing how it acts on morphisms.  The
functoriality and other useful properties of the construction below
are in Proposition 1.1 of \citep{ES} .

Given a map $E(q) \to E(q-p)$ we need to construct a map $\bigwedge^q
\mathcal{T} \to \bigwedge^{q-p} \mathcal{T}$. Any map $E(q) \to E(q-p)$ is
defined by a single element $a \in \bigwedge^p V$.  This also
defines a map $\bigwedge^p W \to \mathbb{C}$.  As $\mathcal{T}$ is a
subbundle of $W \otimes \mathcal{O}_{\mathbf{G_{n+1}}}$, there is an
induced map $a \ :\ \bigwedge^p \mathcal{T} \to
\mathcal{O}_{\mathbf{G_{n+1}}}$. Finally, to construct the map
$U_{n+1}(a) \ : \ \bigwedge^q \mathcal{T} \to \bigwedge^{q-p}
\mathcal{T}$, start with the standard diagonal map $\Delta \ : \
\bigwedge^q \mathcal{T} \to \bigwedge^{q-p} \mathcal{T} \otimes
\bigwedge^{p} \mathcal{T}$ and compose with the map $1 \otimes a$.

We will need to use a more explicit description of the map, in terms
of our chosen bases. Recall that a fiber of $\mathcal{T}$ is a
subspace $W_f = \mathbb{C}\{f_1, \dots, f_{n+1} \}$.  We can write the
$f_i$ as:

$$f_i = \sum_{j=1}^N C_{ij} X_j. $$ 

The coefficients form a $(n+1) \times N$ matrix $C$.  Given ordered
subsets $I = \{i_1, \dots, i_p \} \subset \{1, \dots, n+1 \}$ and $J =
\{j_1, \dots, j_p \} \subset \{1, \dots, N \}$, of the same size $p$,
let $C_{I,J}$ denote the determinant of the submatrix of $C$ with rows
from $I$ and columns from $J$. We will also use the notation $f_I =
(-1)^I \bigwedge_{i \in I} f_i$ and $e_J = \bigwedge_{j
\in J} e_j$. Note the sign factor added to the $f$ part only in order
to simplify the signs in the next proposition:

\begin{lemma} 
\label{l:map}

Let $J \subset \{1, \dots, N \}$ with $|J| = p$. We view $e_J$ as a
map from $E(q)$ to $E(q-p)$. In that case for any $I \subset \{1,
\dots, n+1\}$ with $|I| = q$:
 
$$(U_{n+1}(e_J))(f_I) = \sum_{I_1 \subset I, \ |I_1| = p} C_{I_1, J}
f_{I \setminus I_1} $$

\end{lemma}

\begin{proof} This is a direct translation of the above description
applied to our particular choice of bases. The diagonal map splits
up $f_I$ into a sum of pieces corresponding to a choice of $I_1$ and
its complement. The action of $e_J$ on the piece corresponding to
$I_1$ is exactly the determinant of the specified minor. The only thing
to check is that the sign works out.  
\end{proof}

\section{Toric Tate Resolution}

\label{s:torictate}
We return to the case in question, where $X_A$ is a toric surface with
corresponding polytope $Q$. As we saw earlier, the sections of the
corresponding very ample divisor are just the elements of the vector
space $S_Q$.  Therefore, we will apply the exterior algebra
construction with $W = S_Q$ and $V = S_Q^\ast$. The corresponding
projective space is $\mathbf{P} = \mathbb{P}(W) \cong
\mathbf{P}^{N-1}$, and the exterior algebra is $E = \bigwedge V$.

Any Weil divisor on the toric surface $X_A$ yields a rank one reflexive
sheaf which can be extended to a sheaf on $\mathbf{P}$ under the given
embedding.  We will consider the particular divisor corresponding to
${\rm int}(2Q)$ i.e. $2D - \sum D_i$.  Let $\mathcal{F}$ be the
corresponding sheaf $\mathcal{O}_{X_A}({\rm int}(2Q))$ extended to a
sheaf of $\mathbf{P}$.

\begin{proposition}
\label{p:cohom}
\begin{align}
\label{eqn:H0} H^0 (\mathcal{F}(k)) &\cong S_{{\rm int} ((2+k)Q)} \\ 
\label{eqn:H1} H^1(\mathcal{F}(k)) &\cong 0 \\
\label{eqn:H2} H^2 (\mathcal{F}(k)) &\cong S_{(-2-k)Q}^\ast
\end{align}

\noindent for all $k \in \mathbb{Z}$.
\end{proposition}

\begin{proof}
First of all, since all sheaves are supported on $X_A$, it is
equivalent to compute cohomology on $X_A$.  By construction, $X_A$ is
normal and thus Cohen Macaulay by Hochster's theorem. The dualizing
sheaf is $\mathcal{O}(\omega) = \mathcal{O}(-\sum D_i)$. Also,
twisting by $1$ on $\mathbf{P}$ is the same as twisting by $D$ on
$X_A$. Therefore, $\mathcal{F}(k) = \mathcal{O}((k+2)D - \omega)$.

Now (\ref{eqn:H0}) follows from Proposition \ref{p:sintq}.  For $k > -2$,
$\mathcal{F}(k)$ is an ample divisor minus the canonical
divisor. Therefore, the higher cohomology, $H^1$ and $H^2$ must be
zero by Mustata's vanishing result, \citep[Theorem 2.4 (ii)]{Mus}.

Furthermore, $\mathcal{O}(D)$ is very ample, hence locally free, so
Serre duality tells us $H^i(\mathcal{O}((k+2)D - \omega)) \cong
H^{2-i} (\mathcal{O}((-2-k)D))^\ast$. In particular, applying
Proposition \ref{p:sq} to $i=2$ gives us statement (\ref{eqn:H2}) in the
proposition. For $k \leq -2$, $\mathcal{O}((-2-k)D)$ is generated by
its sections and so all higher cohomology, in particular $H^1$
vanishes, completing the proof of (\ref{eqn:H1}).
\end{proof}

\begin{corollary} The Tate resolution $T(\mathcal{F})$ has terms:

\label{thm:tate}
\begin{align*}
T^e(\mathcal{F}) =& S^\ast_{-eQ} \otimes E(2-e) \quad {\rm for} \ e < -1 \\
T^{-1}(\mathcal{F}) =& S^\ast_Q \otimes E(3) \oplus S_{{\rm int}(Q)} \otimes E(1) \\
T^{0}(\mathcal{F}) =& S^\ast_0 \otimes E(2) \oplus S_{{\rm int}(2Q)} \otimes E(0) \\
T^e(\mathcal{F}) =& S_{{\rm int}(eQ)} \otimes E(-e) \quad {\rm for} \ e >0,
\end{align*}
with maps as follows: 

\begin{tiny}
\begin{diagram}
	\cdots & \rTo^{i_m} & (S_{2Q})^\ast \otimes E(4) & \rTo^{i_m} &
	(S_Q)^\ast \otimes E(3) & \rTo^{i_m} & (S_0)^\ast \otimes E(2)
	& \rTo & 0 & \\
	&& \oplus & \rdTo^{\Delta_{2Q}} & \oplus &
	\rdTo^{\Delta_Q}& \oplus & \rdTo^{\Delta_0} & \oplus &\\
	&& 0 & \rTo & S_{{\rm int}(Q)} \otimes E(1) & \rTo^{\wedge m} &
	S_{{\rm int}(2Q)} \otimes E & \rTo^{\wedge m} & S_{{\rm int}(3Q)} \otimes E(-1) & \rTo^{\wedge m} & \cdots. \\
\end{diagram}
\end{tiny}

The horizontal maps $\wedge m$ and $i_m$ are all multiplication by the
element $m = \sum y^\alpha \otimes e_\alpha$ where $\alpha$ ranges
over the lattice points in $Q$, and $e_\alpha$ is the corresponding
dual vector in $E$. 

\end{corollary}

\begin{proof}

We simply plug in our known cohomology from \ref{p:cohom} into
(\ref{e:tate}) to obtain the terms. The horizontal maps are indeed
multiplication by $m$, as per our discussion in the previous section,
noting only that the Serre duality respects the $S$-module structure
of the cohomology.
\end{proof}

Now we apply the functor $U_3$ to $T(\mathcal{F})$. Once again let
$\mathcal{T}$ denote the tautological bundle on the Grassmannian of
codimension 3 planes in $\mathbb{P}^{N-1}$. Note that $\bigwedge^p
\mathcal{T} = 0$ for $p > 3$ or $p < 0$. Therefore
$U_3(T(\mathcal{F}))$ is  the two term complex below:

\begin{diagram}
	   &      & (S_Q)^\ast \otimes \bigwedge^3 \mathcal{T} 
                & \rTo^{\widehat{i_m}} 
                & (S_0)^\ast \otimes \bigwedge^2 \mathcal{T}
	          &      & \\
	 0 & \rTo & \oplus
                & \rdTo^{\widehat{\Delta}_Q}
                & \oplus &
                & \rTo   & 0.\\  
	   &      & S_{{\rm int}(Q)} \otimes \bigwedge^1 \mathcal{T}
                & \rTo^{\widehat{\wedge m}} 
                & S_{{\rm int}(2Q)}  \otimes \bigwedge^0 \mathcal{T}&\\
\end{diagram}

Since $\mathcal{F}$ is of rank 1, the resultant is up to a constant
 the determinant of the matrix of the nontrivial map $(\widehat{i_m} +
 \widehat{\Delta}_Q) \oplus \widehat{\wedge m}$. However, we can of
 course normalize the maps in the Tate resolution so that we have the
 resultant up to sign. From here on we assume that such a
 normalization has been made.

All that is left to do is describe the maps $\widehat{\wedge m},
\widehat{\Delta}_Q$, and $\widehat{i_m}$. It is enough to define these
maps on each fiber, that is, for each choice of $(f_1, f_2, f_3)$.

To describe the maps $\widehat{\wedge m}$ and $\widehat{i_m}$ we
introduce the {\em Sylvester map} $\Psi_t \ : \ S_t \otimes
\mathbb{C}^3 \to S_{t+Q}$ which sends $(g_1, g_2, g_3)$ to $f_1g_1 +
f_2g_2 + f_3g_3$.

\begin{proposition}
\label{p:sylvmaps}
 The map $\widehat{\wedge m}$ is $\Psi_{{\rm int}(Q)}$, and the map
$\widehat{i_m}$ is $(\Psi_0)^\ast$ on each fiber over the Grassmannian.
\end{proposition}

\begin{proof}
First consider $\wedge m$. We pass to $\bigwedge^1 \mathcal{T}$, which
has a basis at each fiber indexed by $f_1, f_2, f_3$. By Lemma
\ref{l:map}, on the factor corresponding to $f_i$ we must replace each
$e_{\alpha}$ in $m$ by the corresponding coefficient $C_{i\alpha}$.
So on the factor corresponding to $f_i$, multiplication by $m =
\sum_{\alpha \in A} y^\alpha \otimes e_{\alpha}$ becomes
multiplication by $\sum_{\alpha \in A} C_{i\alpha}y^\alpha = f_i$
. This is exactly the Sylvester map.

On the other hand $i_m$ is the map sending $(y^\alpha)^\ast$ to
$e_\alpha$.  To apply the functor $U_3$ we pick the basis $(f_1 \wedge
f_2 \wedge f_3)$ on $\bigwedge^3 \mathcal{T}$ and $(f_2 \wedge f_3, -f_1
\wedge f_3, f_1 \wedge f_2)$ on $\bigwedge^2 \mathcal{T}$. Another
application of Lemma \ref{l:map} shows that $e_\alpha$ is replaced by
the vector $(C_{1\alpha}, C_{2\alpha}, C_{3\alpha})$ in terms of this
second basis. This is exactly the dual Sylvester map $(\Psi_0)^\ast$.
\end{proof}  

Computing $\widehat{\Delta_Q}$ from $\Delta_Q$ is straightforward.

\begin{proposition}

\label{p:bezmap}
Write 

\[ \Delta_Q((y^\alpha)^\ast) = \sum_\beta \sum_{u,v,w} c_{uvw} (e_u \wedge e_v \wedge e_w) y^{\beta}, \]

then for each fiber $(f_1, f_2, f_3)$ on the Grassmannian:

\[ \widehat{\Delta_Q}((y^\alpha)^\ast) =  \sum_\beta \sum_{u,v,w} c_{uvw} [uvw] y^{\beta}. \]

\end{proposition}

\begin{proof}
Here, both $\bigwedge^3 \mathcal{T}$ and $\bigwedge^0 \mathcal{T}$ are 1
dimensional vector spaces. Lemma \ref{l:map} tells us to replace $e_u \wedge e_v \wedge e_w$ by the determinant of the maximal minor with columns $u, v, w$ of the coefficient matrix of the $f_i$, i.e the bracket $[uvw]$. 
\end{proof}   

Putting it all together we have a proof of Theorem \ref{thm:blmtrx}.

\begin{proof}[Proof of Theorem \ref{thm:blmtrx}]
The Chow form is the determinant, up to sign, of the map $(\widehat{i_m} +
\widehat{\Delta}_Q) \oplus \widehat{\wedge m}$. However, the blocks of
the matrix corresponding to $\wedge m$ and $i_m$ are just Sylvester
maps, by Proposition \ref{p:sylvmaps}, whose matrices are $L$ and $\tilde L$
respectively. The matrix of $\widehat{\Delta_Q}$ has entries which are linear
forms in the bracket variables by Proposition \ref{p:bezmap} above.
\end{proof}

As a corollary we note that the matrix must be square. That is, $3 +
\#{\rm int}(2Q) = 3\cdot\#{\rm int}(Q) + \#Q$. This identity also arises from
the simple fact that the third difference of the quadratic Erhart
polynomial of $Q$ is 0.

All that is left is to prove our formula for $\widehat{\Delta_Q}$ in
Theorem \ref{thm:main}, for which, by the above, we need to prove the
corresponding formula for $\Delta_Q$. It turns out that it is easy to
compute $\Delta_0$, and we can verify a formula for $\Delta_Q$ by
making sure it lifts $\Delta_0$. This is described below.

\section{The Map $\Delta_Q$}

\label{s:proofs}

The map $\Delta_0$ is closely related to the {\em toric Jacobian}
\citep{Cox2}. The toric Jacobian is usually constructed as the
determinant of a matrix of partial derivatives.  Cattani, Cox, and
Dickenstein \citep{CCD} construct a different element, which they call
$\Delta_{\sigma}$, referring to the choice $\sigma$ of a cone in the
fan, which is a constant times the Jacobian modulo the ideal $I = (f_1,
f_2, f_3)$. Moreover, while the Jacobian of three forms supported on
$Q$ has {\em toric residue} \citep{CCD} equal to the normalized area of
$Q$, this new element has residue 1. Therefore, we will call this the
{\em normalized Jacobian} and it is unique modulo $I$.

Let $y_1, y_2$ be edge variables such that the
corresponding edges meet at a vertex $p$. Let $y_3, \dots, y_s$ be the
remaining edge variables of the homogeneous coordinate ring $S$. A
monomial $m$ in $S_Q$ is divisible by $y_i$ if and only if the
corresponding lattice point in $Q$ is not on the corresponding edge.

Therefore, we can define a partition of the monomials in $S_Q$ into
three sets $\mu_1, \mu_2, \mu_3$, where $\mu_1$ is defined
to be the set of all monomials divisible by $y_1$, $\mu_2$ is the set
of monomials divisible by $y_2$ but not divisible by $y_1$, and
$\mu_3$ divisible by $y_3\cdots y_s$ but not by either $y_1$ or $y_2$.

Note that $\mu_1$ corresponds to points not on edge 1, $\mu_2$
is the points on edge 1, but not edge 2, and $\mu_3$ is the unique
point, the vertex $p$, on both edges 1 and 2.

\begin{proposition}  Set $M_i = \sum_{s \in \mu_i} s \otimes e_s \in
S_Q \otimes E$. Define $J_0 = \frac{M_1}{y_1} \wedge \frac{M_2}{y_2}
\wedge \frac{M_3}{y_3\cdots y_s}$, an element of $S_{{\rm int}(3Q)}
\otimes E_{-3}$.  A choice for the map $\Delta_0 \ : \ (S_0)^\ast
\otimes E(2) \to S_{int(3Q)} \otimes E(-1)$ is $1 \otimes 1 \mapsto
-J_0$.
\end{proposition}

The element $J_0$ is chosen so that $U_3(J_0)$ is the normalized toric
Jacobian as constructed in \citep{CCD}.

\begin{proof}
First note that the map $\Delta_0$ is determined by the image of $1
\otimes 1 \in (S_0)^\ast \otimes E_0$. By abuse of notation we denote
$\Delta_0(1 \otimes 1)$ by just $\Delta_0$. By exactness, $\Delta_0$
is in the kernel of $\wedge m$ and not in the image of the previous map
$\wedge m$. Furthermore, $\Delta_0$ is unique with respect to this
property, up to a constant and modulo the image of $\wedge m$. Thus we
need to check that our choice $J_0$ is also in the kernel of the
horizontal map $\wedge m$, but not in the image of the previous map
$\wedge m$. Finally, we argue that if we choose the constant -1, the
determinant of the complex will be exactly the resultant (up to sign).

To start with we notice $m = M_1 + M_2 + M_3$, and so $J_0 \wedge m =
\frac{M_1}{y_1} \wedge \frac{M_2}{y_2} \wedge \frac{M_3}{y_3\cdots
y_s} \wedge (M_1 + M_2 + M_3) = 0$. So $J_0$ is indeed in the kernel
of $\wedge m$.

To show that $J_0$ is not in the image of the previous map, we twist
the whole Tate resolution by 1, so that the map $\Delta_0$ goes from
$(S_0)^\ast \otimes E(3)$ to $S_{{\rm int}(3Q)} \otimes E$, and then
apply the functor $U_3$. This also gives a complex whose determinant
is the resultant (Theorem 0.1, in \citep{ES}), in particular it is
exact when the resultant is non-zero. In this situation the image of
the lower map is just the ${\rm int}(3Q)$ graded piece of the ideal $I
= (f_1, f_2, f_3)$, and the normalized toric Jacobian is known to be a
nonzero element modulo this ideal(see \citep{CCD, Cox2}). Therefore,
$J_0$, which specializes to the Jacobian, cannot be in the image of
the map $\wedge m$.

Finally, the specialized complex above, with the normalized toric
Jacobian as the diagonal map, appears in \citep{DE} where the authors
show that the determinant of the complex is exactly the resultant up to
sign. Therefore, the map $1\otimes 1 \mapsto -J_0$ above is a valid
choice, up to sign, for the map $\Delta_0$ in Theorem \ref{thm:tate}.
\end{proof}

Now let's take the degree $-3$ part of the Tate resolution to get:

\begin{tiny}
$$
\begin{diagram}
	0 & \rTo &
	(S_Q)^\ast & \rTo^{i_m} & (S_0)^\ast \otimes \bigwedge^1 V
	& \rTo & 0 \\
	\oplus &  & \oplus &
	\rdTo^{\Delta_Q}& \oplus & \rdTo^{\Delta_0} & \oplus\\
	 0 & \rTo & S_{{\rm int}(Q)} \otimes \bigwedge^2 V & \rTo^{\wedge m} &
	S_{{\rm int}(2Q)} \otimes \bigwedge^3 V & \rTo^{\wedge m}
      & S_{{\rm int}(3Q)} \otimes \bigwedge^4 V. \\
\end{diagram}
$$
\end{tiny}

Let $\{n_\alpha\}$ be the basis of $(S_Q)^\ast$ dual to the monomial
basis $\{y^\alpha\}$ of $S_Q$.  The map on the top row sends
$n_\alpha$ to $e_\alpha$. Because these maps form a complex we have
the relation $\Delta_Q(n_\alpha) \wedge m = -\Delta_0(e_\alpha) =
J_0 \wedge e_\alpha$.

The map $\Delta_Q$ is not canonically defined, even after picking
$\Delta_0$.  In fact the next proposition shows that {\em any} map
satisfying the above relation will do.

\begin{proposition}
Define $\Delta_Q(n_\alpha)$ to be {\em any} element $d_\alpha$,
homogeneous of degree -3, such that $d_\alpha \wedge m = J_0
\wedge e_\alpha$. This defines a valid choice for $\Delta_Q$.
\end{proposition}

\begin{proof}
The map $i_m$ in the top row sending $n_\alpha$ to $e_\alpha$ for each
$\alpha \in Q$ is clearly injective (in fact an isomorphism of vector
spaces). We will use this to show that the bottom row is exact at the
term $S_{{\rm int}(2Q)} \otimes \bigwedge^3 V$. So pick an element $k$
in the kernel of $\wedge m \ : S_{{\rm int}(2Q)} \otimes \bigwedge^3 V
\to S_{{\rm int}(3Q)} \otimes \bigwedge^4 V$. Now $(0,k)$ is in the
kernel of the whole complex. Therefore, by exactness there exists an
element $(a, b) \in (S_Q)^\ast \oplus (S_{{\rm int}(Q)} \otimes
\bigwedge^2 V)$ mapping on to it. But now $i_m(a) = 0$, therefore
$a=0$. So $b \wedge m = k$ as desired. 

Now suppose the Tate resolution is fixed with $\Delta_0$ defined as in
Proposition 6.1. Let $\tilde{\Delta}_Q$ be any map satisfying the
above relation. Therefore, for any $n_{\alpha}$, $\Delta_Q(n_{\alpha})
\wedge m = - {\Delta}_0(e_{\alpha}) = \tilde{\Delta}_Q(n_{\alpha})
\wedge m $. So, $\Delta_Q$ and $\tilde{\Delta}_Q$ differ by an
element of the kernel of $\wedge m$.  By the argument in the previous
paragraph, this is the same as differing by an element of the image of
the previous $\wedge m$.  Therefore, replacing $\Delta_Q$ by
$\tilde{\Delta}_Q$ does not change exactness at this step of the Tate
resolution. As the Tate resolution is a minimal free resolution, this
new choice can always be extended ad infinitum, and so
$\tilde{\Delta}_Q$ is itself a valid map.
\end{proof}

So we need only find for every lattice point $\alpha$ in $Q$, an
element $d_\alpha$ such that $d_\alpha \wedge m = J_0 \wedge
e_\alpha$. In \citep{Khe} it was shown how to reduce this to a problem
in linear algebra. In this paper, we show instead that the explicit,
combinatorial formula from Theorem \ref{thm:main} does the trick. We
restate Theorem \ref{thm:main} below using the language of exterior
algebras developed above. Recall the definitions of the sets $R_i$
from \ref{eqn:Rsets}. The fan has possibly been refined as described
earlier to guarantee that $R_3$ is non-empty.

\begin{theorem}
\label{thm:main2}
 The map $\Delta_Q \ : \ (S_Q)^\ast \otimes E \to S_{{\rm int}(2Q)} \otimes
 E(-3)$ can be defined as follows:

 \[ \Delta_Q(n_\alpha) = \sum_{(u,v,w) \in F_{\alpha} \subset A^3}
 y^{u+v+w - \alpha - \omega_0} \otimes e_u \wedge e_v \wedge e_w. \]

\noindent Here $\omega_0 = (1, 1, \dots, 1)$, and $F_\alpha$ is the set of all
triples $(u,v,w) \in A^3$ satisfying the Boolean combination of inequalities in (\ref{eqn:ineq})

\end{theorem}

The next lemma will rewrite $J_0 \wedge e_{\alpha}$ in a form more
convenient for our purposes.

\begin{lemma}
\label{l:del0alpha} 
\[ J_0 \wedge e_\alpha = \sum_{t,u,v,w}  y^{t+u+v+w - \alpha - \omega_0} \otimes e_u \wedge e_v \wedge e_w \wedge e_t, \]

\noindent where $t, u, v, w$ satisfy:

\begin{align}
\label{ineq1} \forall i \in R_1 \quad t_i + u_i + v_i + w_i >& \alpha_i\\
\label{ineq2} \exists i \in R_1 \quad t_i + v_i + w_i \leq& \alpha_i \\
\label{ineq3} \forall j \in R_2 \quad t_j + v_j + w_j >& \alpha_j \\
\label{ineq4} \exists j \in R_2 \quad t_j + w_j \leq& \alpha_j \\
\label{ineq5} \forall k \in R_3 \quad t_k + w_k >& \alpha_k \\
\label{ineq6} \exists k \in R_3 \quad t_k \leq & \alpha_k.
\end{align}

\end{lemma}

\begin{proof}

First note that if $\exists k \in R_3$ such that $w_k \leq  \alpha_k$, then
both $e_u \wedge e_v \wedge e_w \wedge e_t$ and  $e_u \wedge e_v
\wedge e_t \wedge e_w$, with the same power of $y$, appear in the sum and 
cancel out. So condition (\ref{ineq5}) can be replaced by the stronger condition

\begin{equation}
\forall k \in R_3 \quad w_k > \alpha_k.  \tag{\ref{ineq5}'}
\end{equation}

We will show that every term in $J_0 \wedge e_\alpha$ satisfies these
conditions, and conversely every tuple $(t,u,v,w)$ satisfying the
conditions corresponds to a term in $\Delta_0 \wedge e_\alpha$.

The element $J_0$ can be rewritten as $y^{u+v+w - \omega_0} \otimes
\sum e_u \wedge e_v \wedge e_w $ where $u_1 > 0$, $v_1 = 0$ but $v_2 >
0$, and $w_1 = w_2 = 0$. Wedge this with $e_{\alpha}$, and we show
that the terms $e_u \wedge e_v \wedge e_w \wedge e_{\alpha}$ all
appear on the right hand side. So choose $t = \alpha$ then $t_1 + v_1
+ w_1 = \alpha_1$, $t_2 + w_2 = \alpha_2$ and $t_k = \alpha_k$ for all
k, thus conditions (\ref{ineq2}), (\ref{ineq4}), and (\ref{ineq6}) are
satisfied. On the other hand, $w_i > 0$ for all $i \neq 1,2$. This,
combined with $v_2 > 0$ implies condition (\ref{ineq3}), while $u_1 >
0$ implies condition (\ref{ineq1}). Now, the set $R_3$ is constructed
so that $w$, the vertex where edges 1 and 2 meet, satisfies condition
(\ref{ineq5}') for all $\alpha$ {\em except} when $\alpha = w$, in
which case $J_0 \wedge e_\alpha = 0$. Thus all the terms in $J_0
\wedge e_\alpha$ appear in the desired sum.

Conversely, pick any tuple $(t,u,v,w)$ satisfying (\ref{ineq1}),
(\ref{ineq2}), (\ref{ineq3}), (\ref{ineq4}), (\ref{ineq6}) , and the
modified (\ref{ineq5}'). Define $\gamma = \alpha - t$. So, in our
notation, $\alpha_i - t_i = \langle \eta_i, \gamma \rangle$.

By conditions (\ref{ineq2}), (\ref{ineq4}), (\ref{ineq6}) there exists
$i_0$, $j_0$, $k_0$ in $R_1$, $R_2$, $R_3$ respectively such that
$\langle \eta_{i_0}, \gamma \rangle \geq 0$, $\langle \eta_{j_0},
\gamma \rangle \geq 0$ and $\langle \eta_{k_0}, \gamma \rangle \geq
0$. Since the region $R_3$ is between $R_1$ and $R_2$, we must either
have $\eta_{k_0}$ a positive linear combination of $\eta_{i_0}$ and
$\eta_{j_0}$, or $\langle \eta_{i_0}, \gamma \rangle = \langle
\eta_{j_0}, \gamma \rangle = 0$.

However, we also have $w_{i_0} \leq \alpha_{i_0}$ and $w_{j_0} \leq
\alpha_{j_0}$, but $w_{k_0} > \alpha_{k_0}$, which rules out the first
case.  Thus $t_{i_0} = \alpha_{i_0}$ and $t_{j_0} = \alpha_{j_0}$. By
conditions (\ref{ineq2}) and (\ref{ineq4}) we must have $w_{i_0} =
w_{j_0} = 0$. This is possible only if the facets corresponding to
$\eta_{i_0}$ and $\eta_{j_0}$ meet at a vertex. The only vertex where
the sets $R_1$ and $R_2$ meet is the vertex $p$ when $w_1 = w_2 =
0$. But now, $\gamma$ must be 0, since $\eta_1$ and $\eta_2$ are
linearly independent. Thus $t = \alpha$. So, by condition
(\ref{ineq2}), $v_1 = 0$, by condition (\ref{ineq3}) $v_2 > 0$, and by
condition (\ref{ineq1}), $u_1 > 0$.  Hence, every term in the right
hand sum also appears in $J_0 \wedge e_\alpha$.
\end{proof}

\begin{proof}[Proof of Theorem \ref{thm:main2}] 
We must show that if $\Delta_Q$ is defined as above, then
$\Delta_Q(n_\alpha) \wedge m = J_0 \wedge e_\alpha$.
The left hand side is the sum

\[ \sum_{(u,v,w, t)} y^{u+v+w+t - \alpha - \omega_0} \otimes e_u \wedge e_v \wedge e_w \wedge e_t, \]

\noindent where $(u,v,w)$ satisfy (\ref{eqn:ineq}) and $t$ is unconstrained.

On the other hand, by Lemma \ref{l:del0alpha}, the right hand side
is

\[\sum_{t,u,v,w} y^{t+u+v+w - \alpha - \omega_0} \otimes e_u \wedge e_v \wedge e_w \wedge e_t,  \]

\noindent where $(u,v,w,t)$ satisfy the inequalities
(\ref{ineq1})-(\ref{ineq6}).

So, it is enough to show for any fixed 4 tuple $(u,v,w,t)$ the sum of
all signed permutations satisfying (\ref{eqn:ineq}), is equal to the sum of all
signed permutations satisfying (\ref{ineq1})-(\ref{ineq6}).

We consider the poset corresponding to the power set of $P =
\{u,v,w,t\}$. This is a four-dimensional cube whose vertices are the
16 subsets of $P$, and two subsets $p$ and $q$ are connected by a
directed edge from $p$ to $q$ if $p$ is the union of $q$ with a single
element of $P$.  A maximum oriented path (of length 5) in this poset
corresponds to a permutation of $(u,v,w,t)$. Given a permutation
$(u,v, w, t)$, the path starts at $\emptyset$, has first vertex
$\{t\}$, second vertex $\{w, t\}$ and so on.  Define the {\em sign} of
this path to be the sign of the corresponding permutation. We will
consider formal sums of signed paths, remembering that if the same
path occurs twice in the sum with opposite signs, then the
contribution from that path is 0.

Let $A_i$ be a condition on a vertex $p$ which evaluates to true if
$\sum_{v \in p} v_k > \alpha_k$ holds for all indices $k$ in $R_i$.
Note that if $p$ satisfies $A_i$ and $q \supset p$ then $q$ satisfies
$A_i$. Label a vertex $B_i$ if it satisfies condition $A_i, \dots,
A_3$ but fails to satisfy conditions $A_1, \dots, A_{i-1}$. With this
notation the permutations $(u,v,w,t)$ satisfying (\ref{eqn:ineq}) are
oriented paths through the cube labeled $(B_4, B_3, B_2, B_1,
B_1)$. The permutations, this time ordered $(t, u, v, w)$, satisfying
(\ref{ineq1})-(\ref{ineq6}) are paths of the form $(B_4, B_4, B_3,
B_2, B_1)$. Note that this introduces a sign of $(-1)^3$ into our
formula.

So, to complete the proof it is enough to show the following lemma
that was proved by David Speyer in a personal communication.
\end{proof}

\begin{lemma}
\label{comblemma}
The sum of oriented paths in the cube of the form $(B_4, \dots, B_i,\\
B_i, B_{i-1}, \dots, B_1)$ is $(-1)^{i-1}$ times the sum of paths of
the form  $(B_4, B_3, B_2, \\ B_1, B_1)$.
\end{lemma}

In particular when $i=4$ we have our desired result. 

\begin{proof} 

By induction it is enough to show that the sum of paths of the form
$(B_4, \dots, B_i, B_i, B_{i-1}, \dots, B_1)$ is negative the sum of
paths of the form $(B_4, \dots, \\  B_{i-1}, B_{i-1}, \dots, B_1)$.
Let $S_1$ denote the first sum and $S_2$ the second.

For the moment, consider any two vertices $p$ and $q$ of the cube,
labeled $B_i$ and $B_{i-1}$ respectively, joined by an oriented path
of length 2.  There are exactly two such paths passing through
intermediate vertices $a$ and $b$ respectively. As $a$ contains $p$
and is contained in $q$, by the definition of the labels $a$ satisfies
$A_i, \dots, A_3$ but fails to satisfy $A_1, \dots, A_{i-2}$.  If $a$
obeys $A_{i-1}$ then it has label $B_{i-1}$, otherwise it has label
$B_i$. The case for $b$ is identical.

Returning to the claim consider two disjoint paths of vertices $v_4,
 \dots, v_i$ and $v_{i-1}, \dots v_1$ where $v_j$ has label $B_j$ and
 it is possible to join these paths by adding a single vertex between
 them.  As above, there are two possibilities for this new vertex, $a$
 and $b$, each of which has label $B_i$ or $B_{i-1}$.  The
 permutations associated to the two ways of completing the path differ
 by a single exchange, hence have opposite signs. If $a$ and $b$ have the
same label they cancel in the sum $S_1$ or $S_2$. If they have opposite
labels than one contributes positively to one of the sums, and the
other contributes negatively to the other sum. Therefore, the
two sums are negative of each other.
\end{proof}

\section{Future Work}

\label{s:general}

This paper is, in the author's opinion, just the tip of the iceberg in
the application of exterior algebra methods to sparse resultants. I am
actively working on several more general results and have ideas on
many more.

In this paper we investigated the sheaf $\mathcal{O}({\rm int}(2Q))$
on a toric surface. One of the important properties was the vanishing
of all ``middle'' cohomology. Other sheaves also have this property
and give rise to different formulas for the resultant of a surface. We
can also consider sheaves that do have middle cohomology, although it
seems more difficult to make the maps explicit. In the special case of
products of projective spaces, this is hinted at in Section 6 of the
paper by Dickenstein and Emiris \citep{DiE}.

It is of course of great interest to consider toric varieties of
higher dimension, that is more than 3 equations. I know of a sheaf
giving rise, via the Tate resolution, to a determinantal formula for
the Chow form of any toric threefold. The sticking point is finding an
explicit formula, analogous to Theorem \ref{thm:main}. Hopefully, this
will be worked out in a future publication. 

For four dimensions or higher, it appears the best we can hope for is
matrices whose determinant is a nontrivial multiple of the resultant.
In this situation it should be possible to identify the extraneous
factor with a minor of the matrix. See \citep{DD, D}. 

An important generalization would be to {\em mixed} resultants,
i.e. equations with different supports. Tate resolutions do not
obviously apply, but there may be an appropriate extension.

Finally, returning to the specific formula presented here, there are
several places where choice is involved. An interesting question would
be to classify all possible formulas, for all the different choices.
Another issue is to investigate the efficiency, both in theory and for
an implementation. It may be possible to speed up the computation of
the B\'ezout map $\Delta_Q$.

\bibliographystyle{plain}
\bibliography{jsc02}

\end{document}